\documentclass[11pt]{article}
\usepackage{amsmath}
\usepackage{amssymb}

\newtheorem{thm}{Theorem}[section]

\newtheorem{lemma}[thm]{Lemma}
\newtheorem{prop}[thm]{Proposition}
\newtheorem{rem}[thm]{Remark}
\newtheorem{cor}[thm]{Corollary}
\usepackage{mathpazo}
\usepackage{color}
\begin{document}

\title{Inverses of gamma functions} 

\author{Henrik L. Pedersen\footnote{Research
    supported by grant 10-083122 from The Danish Council for
    Independent Research $\vert$ Natural Sciences}} 
\date{\today}
\maketitle

\begin{abstract}
Euler's Gamma function $\Gamma$ either increases or decreases on intervals between two consequtive critical points. The inverse of $\Gamma$ on intervals of increase is shown to have an extension to a Pick-function and similar results are given on the intervals of decrease, thereby answering a question by Uchiyama. The corresponding integral representations are described. Similar results are obtained for a class of entire functions of genus 2, and in particular integral representations for the double gamma function and the $G$-function of Barnes are found.
\end{abstract}

\section{Introduction}
Euler's gamma function $\Gamma$ increases on the interval $(\alpha,\infty)$, where $\alpha$ denotes the unique zero of the logarithmic derivative $\psi$ of the gamma function on the positive half line.
In \cite{u} the inverse function of $\Gamma$ defined on $(\Gamma(\alpha),\infty)$ was shown to have an extension to a Pick-function in the cut plane $\mathbb C\setminus (-\infty,\Gamma(\alpha)]$. 
In \cite{p} this result was extended to a certain class of entire functions of genus 1. Furthermore, a version of the result was obtained for a class of entire functions of genus 2, including the double gamma function of Barnes.

The construction in these papers relied among other things on the theory of positive and negative definite kernels and in particular on L\"owner's theorem linking positive definite L\"owner kernels with Pick functions.  This method would not seem to work when extending inverses of the gamma function on the interval $(0,\alpha)$ or for that matter on intervals of the negative line where $\Gamma$ has singularities.

The first goal of this paper is to extend the inverses of the gamma function on the remaining intervals of the real line to holomorphic functions in the upper half plane and describe those as  Pick-functions. This, in particular, answers a question posed in \cite{u} about the properties of the analytic extension of the inverse of $\Gamma$ on $(0,\alpha)$. See Corollary \ref{cor:uchiyama}. We shall base our construction on results relating functions of Laguerre-P\'olya class to conformal mappings of the upper half plane onto certain ``comb'' domains. 

In sections \ref{sec:logGamma}, \ref{sec:construction} and \ref{sec:integral} the case of the Gamma function is described in detail, including a description of the integral representations. See theorems \ref{thm:mainodd} and  \ref{thm:maineven}. The method can also be used to study inverses of Laguerre-P\'olya functions, but we do not aim at a complete description of the properties of inverses of Laguerre-P\'olya functions. In Remark \ref{rem:sin} the situation for inverses to $\sin z$ is sketched briefly.

The second goal of this paper is to investigate inverses for a class of entire functions of genus 2, not in the Laguerre-P\'olya class. This class includes the double gamma function and the $G$-function of Barnes. We show that inverses of these functions initially defined on a positive half-line can be extended to Pick-functions. See Section \ref{sec:genus2}.

A Pick-function is a holomorphic function $p$ defined in the upper half plane $\mathbb C_+$ such that $\Im p(z)\geq 0$ for all $z\in \mathbb C_+$.

An entire function $f$ is said to belong to the
Laguerre-P\'olya class if it is the limit of a sequence of real polynomials with all zeros real. This class can be characterized in terms of the Weierstra\ss\ factorization and also in terms of conformal mappings of the upper half plane onto so-called $\mathcal V$-comb domains. 

Suppose that $m,k\in \mathbb Z\cup \{\pm \infty\}$ and $-\infty\leq m<k\leq \infty$. Furthermore let $\{h_j\}$ be a sequence in $[-\infty,\infty)$. The $\mathcal V$-comb domain corresponding to $m,k,\{h_j\}$ is the domain  of the form 
$$
\mathcal V=\{ z\, |\, m\pi <\Im z<k\pi\} \setminus \cup_{m<j<k}\{x+i\pi j\, |\, x\leq h_j\}.
 $$
 The connection to conformal mapping goes back to papers of Grunsky (\cite{g}), Vinberg (\cite{v}), and MacLane (\cite{m}). See also the preprint \cite{kk}. We shall use the following fact from a recent paper by Eremenko and Yuditskii, see \cite[Theorem 2.3]{ey}: 

If $f$ is an entire function of Laguerre-P\'olya class with Weierstra\ss\ factorization
$$
f(z)=z^qe^{-az^2+bz}\prod_{k=1}^{\infty}(1-z/z_k)e^{z/z_k},
$$
where $a\geq 0$, $b\in \mathbb R$, $q\in \{0,1,\cdots\}$, and $\sum_{k=1}^{\infty}|z_k|^{-2}<\infty$ then $f(z)=e^{\phi(z)}$, where $\phi:\mathbb C_+\to \mathcal V$ is a conformal maping onto a $\mathcal V$-comb.

\section{The conformal mapping $\log \Gamma$}
\label{sec:logGamma}
It is a classical fact that the reciprocal of the gamma function belongs to the Laguerre-P\'olya class and admits the following Weierstra\ss\ factorization
$$
\frac{1}{\Gamma(z)}=ze^{\gamma z}\prod_{k=1}^{\infty}\left(1+z/k\right)e^{-z/k}.
$$
The function
\begin{equation}
\label{eq:logGamma}
\log \Gamma(z)=-\gamma z-\log z-\sum_{k=1}^{\infty}\log \left(1+z/k\right)-z/k
\end{equation}
is holomorphic in the cut plane $\mathbb C\setminus (-\infty,0]$ and extends $\log \Gamma(x)$ for $x>0$. (The logarithms on the right hand side of the relation are the principal logarithm defined in the cut plane.)

The logarithmic derivative $\psi=\Gamma'/\Gamma$ is commonly known as the psi-function and it has the representation
$$
\psi(z)=-\gamma+\sum_{k=0}^{\infty}\left(\frac{1}{k+1}-\frac{1}{k+z}\right).
$$
Clearly, $\psi$ increases strictly from $-\infty$ to $\infty$ on each interval $(-k,-k+1)$, for $k\geq 1$. There is thus a unique point $x_k\in (-k,-k+1)$ such that $\Gamma'(x_k)=0$. When $k$ is even $\Gamma$ attains its (positive) minimum value on $(-k,-k+1)$ at $x_k$ and when $k$ is odd its (negative) maximum at $x_k$. 
We put $x_0=\alpha$, the minimum point of $\Gamma$ on the positive half line.

By the results above $\log \Gamma$ defines a conformal mapping of the upper half plane onto the domain $\mathcal V$ of the form 
\begin{equation}
\label{eq:V}
\mathcal V=\mathbb C\setminus \cup_{k=0}^{\infty} [\log|\Gamma(x_k)|,\infty)\times \{-ik\pi\}.
\end{equation}
We notice that $\log \Gamma$ extends to a continuous function on the closed upper half plane except the points $\{0,-1,-2,\cdots\}$, indeed
$$
\lim_{z\to x,\, z\in \mathbb C_+ }\log \Gamma(z)=\log|\Gamma(x)|-ik\pi
$$
for $x\in (-k,-k+1)$, $k\geq 1$. 

We remark that conformal properties of the function $\log \Gamma$ in right half planes was investigated in \cite{ar}. In a rather old paper by Ginzel (see \cite{gi}) mapping properties of $\Gamma$ were also described. 

\begin{rem}
\label{rem:argGamma}
From \eqref{eq:logGamma} it follows that 
$$
\partial_x\arg \Gamma (x+iy)=\Im \psi(x+iy)=\sum_{k=0}^{\infty}\frac{y}{(k+x)^2+y^2}
$$
so for fixed $y>0$, $x\mapsto \arg \Gamma (x+iy)$ is strictly increasing. Since 
$$
m\partial_x\arg \Gamma (m+iy)\geq  m\sum_{k=0}^{m}\frac{y}{(k+m)^2+y^2}\geq \frac{ym^2}{4m^2+y^2}\geq \frac{y}{4+y^2}
$$
for $m\geq 1$ it follows that $\lim_{x\to \infty}\partial_x\arg \Gamma (x+iy)=\infty$. Furthermore, using the functional equation we see that $\arg \Gamma(-m+iy)\leq -m\pi/2+\arg \Gamma(iy)$ for $m\geq 0$, so  $\lim_{x\to -\infty}\partial_x\arg \Gamma (x+iy)=-\infty$. Hence $x\mapsto \arg \Gamma (x+iy)$  maps $\mathbb R$ onto $\mathbb R$.

Finally, we record another consequence of \eqref{eq:logGamma}: $\arg \Gamma(iy)\to \infty$ as $y\to \infty$. 
\end{rem}

\section{Construction of inverses}
\label{sec:construction}
For an integer $k$ let $S_k$ denote the strip
$$
S_k=\{z\in \mathbb C| \Im z\in (-(k+1)\pi,-k\pi)\}.
$$
We consider 
$$
g_k(z)=(\log \Gamma)^{-1}\left(\log z-i(k+1)\pi\right), \quad z\in \mathbb C_+.
$$
This is a holomorphic function in the upper half plane. Let $\mathcal D_k=g_k(\mathbb C_+ )$ be the conformal image of $\mathbb C_+$. Then, clearly,
$$
\mathcal D_k=\{z\in \mathbb C_+ | \log \Gamma (z)\in S_k\}=\{z\in \mathbb C_+ | -(k+1)\pi<\arg \Gamma (z)<-k\pi\}.
$$
\begin{prop}
\label{prop:pick}
For all $k\in \mathbb Z$, $g_k$ is a Pick-function and $\Gamma(g_k(z))=(-1)^{k+1}z$ for $z\in \mathbb C_+ $.
\end{prop}
{\it Proof.} Let $\mathcal V$ be the domain defined in \eqref{eq:V}. Since $\log$ maps the upper half plane onto the strip $S_{-1}$, $\log z-i(k+1)\pi \in \mathcal V$ for $z\in \mathbb C_+ $, and hence $(\log \Gamma)^{-1}\left(\log z-i(k+1)\pi\right)\in \mathbb C_+ $. By construction 
$$
\Gamma(g_k(z))=e^{\log \Gamma(g_k(z))}=e^{\log z-i(k+1)\pi}=(-1)^{k+1}z.
$$
This proves the proposition. \hfill $\Box$
\begin{prop}
The restriction of $\Gamma$ to $\mathcal D_k$ is conformal.
\end{prop}
{\it Proof.} If $\Gamma(z_1)=\Gamma(z_2)$ for $z_1,z_2\in \mathcal D_k$ then $\log|\Gamma(z_1)|=\log|\Gamma(z_1)|$ and $\arg \Gamma(z_1)=\arg \Gamma(z_2)+2\pi l$, for some $l\in \mathbb Z$. But $\arg \Gamma(z_1)$ and $\arg \Gamma(z_2)$ can only differ at most $\pi$, since both $z_1$ and $z_2$ belong to $\mathcal D_k$. Hence $l=0$ and thus $\log \Gamma(z_1)=\log \Gamma(z_2)$. Since $\log \Gamma$ is conformal this yields $z_1=z_2$.\hfill $\Box$

The function $g_{-1}$ maps the upper half plane onto $\mathcal D_{-1}$ and is by reflection in the real line holomorphic in $\mathbb C\setminus \mathbb R$. It has a holomorphic extension to $\mathbb C\setminus (-\infty, \Gamma(x_{0})]$: this can be seen by noting that $g_{-1}$ has a continuous and real valued extension to the interval $(\Gamma(x_{0}), \infty)$ and then conclude using Morera's theorem. One can also use an argument based on the fact that $\Gamma$ has a local inverse at each point in $(\Gamma(x_{0}), \infty)$. See the proof of Proposition \ref{prop:G_k}.

Uchiyama calls the function $g_{-1}$ the principal inverse of the gamma function. From Proposition \ref{prop:pick}
the result of Uchiyama follows:
\begin{cor}
The principal inverse of $\Gamma$ is a Pick-function.
\end{cor}

We shall go through similar constructions on the remaining intervals $(x_{k+1},x_k)$ for $k\geq 0$.

We notice that $\Gamma(x_{2l-1})<0<\Gamma(x_{2l})$ for $l\geq 0$ and define the interval $I_k$ as 
$$
I_k=\left\{ \begin{array}{ll}
[\Gamma(x_k),\Gamma(x_{k+1})], & k\ \text{odd}\\ { }
[-\Gamma(x_{k}),-\Gamma(x_{k+1})], & k\ \text{even.}
\end{array}
\right.
$$
When $k$ is odd, $\Gamma$ increases on each of the intervals $(x_{k+1},-k)$ and $(-k,x_k)$, and $\Gamma((x_{k+1},-k)\cup(-k,x_k))=\mathbb R\setminus I_k$. When $k$ is even, $\Gamma$ decreases on $(x_{k+1},-k)$ and on $(-k,x_k)$, and 
$\Gamma((x_{k+1},-k)\cup(-k,x_k))=\mathbb R\setminus (-I_k)$. 

For $k\geq 0$ $g_k$ is extended by reflection in the real line to a holomorphic function in $\mathbb C\setminus \mathbb R$ and mapping $\mathbb C\setminus \mathbb R$ onto $\mathcal D_k\cup \overline{\mathcal D_k}$. We now show that $g_k$ can be extended across $\mathbb R\setminus I_k$ where $I_k$ is given above. 

\begin{prop}
\label{prop:G_k}
Let $k\geq 0$. The function $g_k$ has a holomorphic extension $G_k$ to $\mathbb C\setminus I_k$ satisfying $G_k(\overline{z})=\overline{G_k(z)}$.
\end{prop}
{\it Proof.} Assume that $k$ is odd. (The argument for $k$ even is similar.) Let $y$ be any real number less than $\Gamma(x_k)$. There exists a unique $x\in (-k,x_k)$ such that $\Gamma(x)=y$. Since $\Gamma'(x)\neq 0$, $\Gamma$ is locally one-to-one and thus there are neighbourhoods $U$ of $x$ and $V$ of $y$ such that for any $w\in V$ there is a unique $z\in U$ satisfying $\Gamma(z)=w$ and the mapping $\Gamma^{-1}:w\mapsto z$ is holomorphic in $V$. Since $\Gamma^{-1}(V\cap \mathbb R)\subseteq \mathbb R$ it follows that  $\Gamma^{-1}(\overline{w})=\overline{\Gamma^{-1}(w)}$ for $w\in V$.

A similar argument shows that for any $y>\Gamma(x_{k+1})$ there exist $x\in (x_{k+1},-k)$ and neighbourhoods of $U$ of $x$ and $V$ of $y$ and a holomorphic function $\Gamma^{-1}: V\to U$ such that $\Gamma^{-1}(\overline{w})=\overline{\Gamma^{-1}(w)}$ for $w\in V$.

Defining $G_k(w)$ as $g_k(w)$ for $w$ in the upper half plane, as $\overline{g_k(\overline{w})}$ for $w$ in the lower half plane and as $\Gamma^{-1}(w)$ for $w\in \mathbb R\setminus I_k$ it follows that $G_k$ is holomorphic in  $\mathbb C\setminus I_k$. \hfill $\Box$
 
\begin{rem}
\label{rem:sin}
The trigonometric functions $\sin$ and $\cos$ are also members of the Laguerre-P\'{o}lya class. It follows that $\log \sin z$ (defined via the infinite product and the principal logarithms) maps $\mathbb C_+$ conformally onto the $\mathcal V$-comb
$$
\mathbb C\setminus \left\{\cup_{k\in\mathbb Z}(-\infty,0]\times \{ik\pi\}\right\}.
$$
Hence, an inverse to $\sin$ can be constructed by $g(z)=(\log \sin)^{-1}(\log z)$ and this inverse maps $\mathbb C_+$ onto the half strip $\{ |\Re z|< \pi/2, \Im z>0\}$. 

The function $g$ can be representated in terms of the inverse to the Joukowski transformation: The Joukowski transformation is defined by $J(w)=(w+1/w)/2$ and it maps $\mathbb C_+\setminus \{|z|\leq 1\}$ conformally onto $\mathbb C_+$. Hence 
$$
z\mapsto i\log(J^{-1}(z))+\pi/2
$$ 
also maps $\mathbb C_+$ conformally onto $\{ |\Re z|< \pi/2, \Im z>0\}$. It is easily seen that the two maps agree on $i\mathbb R_+$ and therefore 
$g(z)=i\log(J^{-1}(z))+\pi/2$ is an inverse to $\sin$ when extended holomorphically from $(-1,1)$ to $(-\pi/2,\pi/2)$. (This formula is not new, see e.g.\ \cite[1.622]{gr}.)
\end{rem}

\section{Integral representation}
\label{sec:integral}
It is well known that any Pick function $p$ has an integral representation of the form 
$$
p(z)=az+b+\int_{-\infty}^{\infty}\left(\frac{1}{t-z}-\frac{t}{t^2+1}\right)\, d\mu(t),
$$
where $a\geq 0$, $b\in \mathbb R$, and $\mu$ is a positive measure on $\mathbb R$ such that 
$$
\int_{-\infty}^{\infty}\frac{d\mu(t)}{t^2+1}<\infty.
$$
Furthermore, $a=\lim_{y\to \infty}p(iy)/y$, $b=\Re p(i)$, and 
$$
\mu =\frac{1}{\pi}\lim_{y\to 0_+}\Im p(x+iy)dx
$$
in the vague topology.

For the function $p=g_{-1}$, the corresponding measure $\mu_{-1}$ is concentrated on $(-\infty, \Gamma(x_{0})]$ and for $k\geq 0$ the measure $\mu_k$ corresponding to $g_{k}$ is supported on $I_k$. We investigate properties of these  measures in more detail.

\begin{lemma} 
\label{lemma:mu_k}
The measure $\mu_k$ in the representation 
$$
g_k(z)=a_kz+b_k+\int_{-\infty}^{\infty}\left(\frac{1}{t-z}-\frac{t}{t^2+1}\right)\, d\mu_k(t),
$$
is supported on $I_k$ and on $I_k\setminus \{0\}$ it has a $C^{\infty}$-density $d_k$ given as 
\begin{equation}
\label{def:d_k}
d_k(t)\equiv \frac{1}{\pi}\Im g_k(t+0i), \quad t\in I_k\setminus\{0\}.
\end{equation}
\end{lemma}
{\it Proof.} We shall give the proof in the situation where $k\geq 0$ is odd (in which case $\Gamma(x_k)<0<\Gamma(x_{k+1})$).
It is convenient to consider two holomorphic extensions of the function $g_k$. The first one is defined by 
$$
h^+_k(z)=\left(\log \Gamma\right)^{-1}(\log z-i\pi (k+1)).
$$
This is a holomorphic extension of $g_k$ to $\mathbb C\setminus \left((-\infty, 0]\cup [\Gamma(x_{k+1}),\infty)\right)$. The second one is defined as
$$
h^-_k(z)=\left(\log \Gamma\right)^{-1}(\log(-z)-i\pi k),
$$
and since $\log (-z)=\log (z)-i\pi$ for $z\in \mathbb C_+ $,  $h^-_k$ is a holomorphic extension of $g_k$ to $\mathbb C\setminus \left((-\infty, \Gamma(x_{k})]\cup [0,\infty)\right)$.

Let $\phi$ be a continuous function of compact support in $I_k\setminus \{0\}$. Then, using dominated convergence, 
\begin{align*}
\lim_{y\to 0_+}\frac{1}{\pi}&\int_{-\infty}^{\infty}\phi(x)\Im g_k(x+iy)dx\\
=&\int_{\Gamma(x_{k})}^0\phi(x) h^-_k(x)dx+\int_0^{\Gamma(x_{k+1})}\phi(x) h^+_k(x)dx. 
\end{align*}
Therefore $\mu_k$ has a $C^{\infty}$-density $d_k$ wrt.\ Lebesgue measure  given as $d_k(t)=\Im g_k(t+0i)/\pi$ for   $t\in I_k\setminus\{0\}$.  (The density is of course integrable wrt.\ Lebesgue measure on $I_k$.) \hfill $\square$
\begin{prop}
\label{prop:d_k-increasing}
For odd $k\geq 1$, the function $d_k$ increases from $0$ to $\infty$ on $I_k\cap (-\infty,0)$ and decreases from $\infty$ to 0 on  $I_k\cap (0,\infty)$. 
\end{prop}
{\it Proof.} Since $\Gamma(g_k(z))=z$ we infer that 
$g_k'(z)=1/(\psi(g_k(z))z)$,
and since $\psi$ is a Pick-function, $\partial_t\Im g_k(t)$ is negative for $t\in I_k\cap (0,\infty)$. 

Let $y_0$ be a given positive number. From Remark \ref{rem:argGamma} it follows that we may choose $x_0\in \mathbb R$ such that $\arg \Gamma(x_0+iy_0)=-\pi(k+1)$. Hence $\log \Gamma(x_0+iy_0)=\log |\Gamma(x_0+iy_0)|-i\pi(k+1)\in \mathcal V$, where $\mathcal V$ is given in \eqref{eq:V}. Since $k+1$ is even it follows that $\Gamma(x_0+iy_0)\in (0,\Gamma(x_{k+1}))$ and also  $\Im g_k(\Gamma(x_0+iy_0))=y_0$.  

Similar arguments show that $\Im g_k$ is increasing from 0 to $\infty$ to the left of the origin.\hfill $\Box$

\begin{thm}
\label{thm:mainodd}
For odd $k\geq 1$ the function $g_k$ has the integral representation
$$
g_k(z)=\int_{\Gamma(x_k)}^{\Gamma(x_{k+1})}\frac{d_k(t)}{t-z}dt-k, \quad z\in \mathbb C\setminus I_k,
$$
where $d_k$ is given in \eqref{def:d_k}. 
\end{thm}
{\it Proof.} From Lemma \ref{lemma:mu_k} we have  
$$
g_k(z)=a_kz+b_k+\int_{\Gamma(x_k)}^{\Gamma(x_{k+1})}\frac{d_k(t)}{t-z}dt-\frac{c_k}{z}, \quad z\in \mathbb C\setminus I_k,
$$
where $a_k\geq 0$, $b_k\in \mathbb R$, and $c_k\geq 0$. Thus $g_k$ is strictly increasing on $\mathbb R\setminus I_k$. Since $\lim_{x\to \infty}g_k(x)=-k$ we obtain $a_k=\lim_{x\to \infty}g_k(x)/x=0$ and then $b_k=-k$. 

The possible point mass $c_k$ at the origin can be computed via 
$$
c_k=\lim_{y\to 0+}y\Im g_k(iy).
$$
We claim that $c_k=0$. Suppose that $c_k>0$. Then $\Im g_k(iy)\to \infty$ as $y\to 0_+$. From this we obtain that $\Re g_k(iy)\leq 0$ for all $y$ sufficiently close to $0$: if for arbitrary small and positive $y$ we would have $\Re g_k(iy)>0$ then 
$$
\arg \Gamma(\Re g_k(iy)+i\Im g_k(iy))>\arg \Gamma(i\Im g_k(iy))
$$
since $\arg \Gamma$ increases on horizontal lines (see Remark \ref{rem:argGamma}). By the same remark, $ \arg \Gamma(iY)\to \infty$ as $Y\to \infty$ and this contradicts the fact that 
$$
-(k+1)\pi<\arg \Gamma(g_k(iy))<-k\pi.
$$ 
Hence we are in the situation that $\Im g_k(iy)\to \infty$ and $\Re g_k(iy)\leq 0$. The next step is to use Stirling's formula
$$
\log \Gamma(w)=\log \sqrt{2\pi}+(w-1/2)\log w-w+\mu(w),
$$
where
$$
\mu(w)=\frac{1}{2}\int_0^{\infty}\frac{Q(t)}{(w+t)^2}\,dt,
$$
and $Q$ is the 1-periodic function defined by $Q(t)=t-t^2$ for $t\in (0,1)$. It is a fact that $|\mu(w)|\leq \pi/8$ for $w\in \mathbb C_+ \setminus \{\Re w\leq 1, 0<\Im  w\leq 1\}$.

By construction of $g_k$ we have 
$$
\log \Gamma(g_k(z))=\log z-i(k+1)\pi.
$$
The idea is now to plug $z=iy$ into this relation, and take real parts on both sides. It gives us
\begin{align*}
\log y=&\log \sqrt{2\pi}+(\Re g_k(iy)-1/2)\log |g_k(iy)|\\
&-\Im g_k(iy)\arg g_k(iy)-\Re g_k(iy)+\Re \mu(g_k(iy)).
\end{align*}
After multiplication by $y$ and some rearrangements we obtain:
\begin{align*}
y\Im g_k(iy)\arg g_k(iy)=&\, -y\log y+y\log \sqrt{2\pi}-\frac{y}{2}\log |g_k(iy)|+y\Re \mu(g_k(iy))\\
&\, +y\Re g_k(iy)(\log |g_k(iy)|-1).
\end{align*}
Since $|g_k(iy)|\to  \infty$, $\Re g_k(iy)\leq 0$ as $y\to 0_+$ this yields
$$
y\Im g_k(iy)\arg g_k(iy)\leq  |y\log y|+y\log \sqrt{2\pi}+\pi y/8.
$$
Since, furthermore, $\arg (g_k(iy))\in [\pi/2,\pi)$ this gives an upper bound on $ y\Im g_k(iy)$, from which we see that this quantity must tend to zero. This is a contradiction and we have shown that $c_k=0$.
\hfill $\Box$

\begin{rem}
From Theorem \ref{thm:mainodd} it follows that, by monotone convergence, 
$$
\int_{\Gamma(x_k)}^{\Gamma(x_{k+1})}\frac{d_k(t)}{x-t}dt\uparrow \int_{\Gamma(x_k)}^{\Gamma(x_{k+1})}\frac{d_k(t)}{\Gamma(x_{k+1})-t}dt 
$$
as $x\downarrow \Gamma(x_{k+1})$ and at the same time $g_k(x)\downarrow x_{k+1}$. Thus 
$$
x_{k+1}=\int_{\Gamma(x_k)}^{\Gamma(x_{k+1})}\frac{d_k(t)}{t-\Gamma(x_{k+1})}dt -k
$$
Similarly, $$
x_{k}=\int_{\Gamma(x_k)}^{\Gamma(x_{k+1})}\frac{d_k(t)}{t-\Gamma(x_{k})}dt-k.
$$
The integral representation of $g_k$ can thus be written as
\begin{align*}
g_k(z)&=x_{k+1}+(z-\Gamma(x_{k+1}))\int_{\Gamma(x_k)}^{\Gamma(x_{k+1})}\frac{d_k(t)dt}{(t-z)(t-\Gamma(x_{k+1}))}\\
&=x_{k}+(z-\Gamma(x_{k}))\int_{\Gamma(x_k)}^{\Gamma(x_{k+1})}\frac{d_k(t)dt}{(t-z)(t-\Gamma(x_{k}))}.
\end{align*}
\end{rem}

\begin{prop}
\label{prop:end_point}
For odd $k\geq 1$ both $d_k(\Gamma(x_k)+s)$ and $d_k(\Gamma(x_{k+1})-s)$ are $O(s^{1/2})$ as $s\to 0_+$.
\end{prop}
{\it Proof.} Since $\Gamma'(x_{k+1})=0$ and $\Gamma''(x_{k+1})>0$ there exists a holomorphic function $m(z)$ in a neighbourhood of $x_{k+1}$ such that 
$$
\Gamma(z)=\Gamma(x_{k+1})+m(z)^2
$$ 
and $m'(x_{k+1})>0$. The function $m$ is locally one-to-one so there exist neighbourhoods $U$ of $x_{k+1}$ and $V$ of $0$ and a holomorphic function $p:V\to U$ such that 
$p(0)=x_{k+1}$, $p'(0)>0$ and $m(p(w))=w$. This gives
$$
\Gamma(p(w))=\Gamma(x_{k+1})+w^2.
$$ 
Since $p'(0)>0$ we may assume that $V$ is chosen so that the points in $V$ intersected with the first quardrant are mapped by $p$ into the domain $\mathcal D_k$.  By construction of the inverse $g_k$, $\Gamma(g_k(w^2+\Gamma(x_{k+1})))=\Gamma(x_{k+1})+w^2$ for $w\in V$ having positive real part. Since $\Gamma$ is conformal in $\mathcal D_k$ we must have $g_k(w^2+\Gamma(x_{k+1}))=p(w)$ for all $w\in V$ having positive real part. Letting $w\to i\sqrt{t}$ for small $t>0$ we obtain $\pi d_k(-t+\Gamma(x_{k+1}))=p(i\sqrt{t})=O(\sqrt{t})$.
 
The investigation at the other end point $\Gamma(x_k)$ is similar, writing $\Gamma(z)=\Gamma(x_{k})-m(z)^2$, where  
$m$ is chosen so that $m'(x_{k})<0$.\hfill $\Box$

In the rest of this section the situation where $k$ is even is briefly described.
\begin{thm}
\label{thm:maineven}
For even $k\geq 0$ the function $g_k$ has the integral representation
$$
g_k(z)=\int_{-\Gamma(x_{k})}^{-\Gamma(x_{k+1})}\frac{d_k(t)}{t-z}dt-k, \quad z\in \mathbb C\setminus I_k,
$$
where $d_k$ is given as $d_k(t)=\frac{1}{\pi}\Im g_k(t+i0)$.
\end{thm}
By considering the function $e_k(z)=g_k(-z)$ in stead of $g_k$ we find:
\begin{cor}
\label{cor:uchiyama}
For even $k\geq 0$ the inverse of $\Gamma: (x_{k+1},-k)\cup (-k,x_k)\to (-\infty, \Gamma(x_{k+1}))\cup (\Gamma(x_k),\infty)$ has an extension $e_k$ to $\mathbb C\setminus [\Gamma(x_{k+1}),\Gamma(x_k)]$ with the representation 
$$
e_k(z)=\int_{\Gamma(x_{k+1})}^{\Gamma(x_{k})}\frac{d_k(-t)}{z-t}dt-k, \quad z\in \mathbb C\setminus [\Gamma(x_{k+1}),\Gamma(x_k)].
$$
(The imaginary part of $e_k$ is negative in the upper half plane.)
\end{cor}
\begin{rem}
This corollary describes for $k=0$ the extension $e_0$ of the inverse of the restriction of $\Gamma$ to $(0,\Gamma(x_0))$ and thus answers the question posed in \cite{u}.
\end{rem}

\section{Entire functions of genus 2}
\label{sec:genus2}
Let $f$ be an entire function having the representation 
\begin{equation}
\label{eq:G-class}
f(z)=z^re^{az^2+bz}\prod_{k=1}^{\infty}\left(1+\frac{z}{\lambda_k}\right)e^{-z/\lambda_k+z^2/(2\lambda_k^2)},
\end{equation}
where $r$ is a non negative integer, $a$ and $b$ are real numbers and the sequence $\{\lambda_k\}$ satisfies $0<\lambda_1\leq \lambda_2\leq \cdots$, and
$$
\sum_{k=1}^{\infty}\frac{1}{\lambda_k^2}=\infty,\quad \sum_{k=1}^{\infty}\frac{1}{\lambda_k^3}<\infty.
$$
This is expressed as the sequence $\{\lambda_k\}$ be of rank 2. The holomorphic function $\log f$ is defined in the cut plane as 
\begin{equation}
\label{eq:logf}
\log f(z)=r\log z +az^2+bz+\sum_{k=1}^{\infty}\log \left(1+\frac{z}{\lambda_k}\right)-\frac{z}{\lambda_k}+\frac{z^2}{2\lambda_k^2}.
\end{equation}

In \cite[Lemma 3.1]{p} it was shown that if $r>0$ then there exists $u>0$ with the property that $(\log f)''(u)=0$, $\log f(x)$ is concave for $0<x<u$ and convex for $x>u$. From this a mapping property can be deduced:
\begin{lemma}
\label{lemma:conformal}
Let $f$ satisfy \eqref{eq:G-class} and suppose that $r>0$. Let $u$ denote the unique point in $(0,\infty)$ such that $(\log f)''(u)=0$. Then:
\begin{enumerate}
\item[(i)] $(\log f)'$ is univalent in
$\{\Re z>u\}$ and maps  $\{\Re z>u\}\cap \mathbb C_+$ into $\mathbb C_+$;
\item[(ii)] $\log f$ is univalent in $\{\Re z>u\}\cap \mathbb C_+$ and in
  $\{\Re z>u\}\cap \mathbb C_-$.
\end{enumerate}
\end{lemma}
{\it Proof.} 
From the representation \eqref{eq:logf} of
$\log f$ it follows that 
\begin{align}
\Im ( (\log f)'(&x+iy))\nonumber\\
=y &\left(-\frac{r}{x^2+y^2}+2a+\sum_{k=1}^{\infty}\left(\frac{1}{\lambda_k^2}-\frac{1}{(x+\lambda_k)^2+y^2}\right)\right),
\label{eq:imlogf'}
\end{align}
and also that
\begin{equation}
\label{eq:logf''}
(\log
f)''(x+iy)
=-\frac{r}{(x+iy)^2}+2a+\sum_{k=1}^{\infty}\left(\frac{1}{\lambda_k^2}-\frac{1}{(x+iy+\lambda_k)^2}\right).
\end{equation}
Therefore, assuming $y\geq 0$, 
\begin{align*}
\Im \left( (\log
f)'(x+iy)\right)
&\geq
y\left(-\frac{r}{x^2}+2a+\sum_{k=1}^{\infty}\left(\frac{1}{\lambda_k^2}-\frac{1}{(x+\lambda_k)^2}\right)\right)\\
&= y\, (\log f)''(x).
\end{align*}
Thus convexity of $\log f$ on $(u,\infty)$ yields that $(\log f)'$
maps $\{\Re z>u\}\cap \mathbb C_+$ into $\mathbb C_+$ and hence also
$\{\Re z>u\}\cap \mathbb C_-$ into $\mathbb C_-$. 

From \eqref{eq:logf''} it follows that 
$(\log f)''$ maps $\{\Re z>0\} \cap \mathbb C_+$
into $\mathbb C_+$. For  $z_1,z_2\in \{\Re z>0\} \cap \mathbb
H_+$, 
$$
(\log f)'(z_2)-(\log f)'(z_1)=(z_2-z_1)\int_0^1(\log f)''(z_1+t(z_2-z_1))\, dt
$$
and hence $(\log f)'$ is univalent in $\{\Re z>0\} \cap \mathbb C_+$. A
similar argument shows that $(\log f)'$ is univalent in $\{\Re z>0\} \cap
\mathbb C_-$. Furthermore, $(\log f)'$ is
increasing on $(u,\infty)$ and all together it means that $(\log f)'$
is univalent in $\{\Re z>0\}$. 

Finally, for $z_1,z_2\in \{\Re z>u\} \cap \mathbb C_+$,
\begin{equation}
\label{eq:logf-int}
\log f(z_2)-\log f(z_1)=(z_2-z_1)\int_0^1(\log f)'(z_1+t(z_2-z_1))\, dt
\end{equation}
and as before this yields that $\log f$ is univalent in $\{\Re z>u\} \cap \mathbb C_+$ and hence also in $\{\Re z>u\} \cap \mathbb C_-$. \hfill $\square$

A function $f$ with the representation \eqref{eq:G-class} either increases on the positive real line or decreases on some open bounded interval $(\alpha,\beta)$, where $0<\alpha$ and increases on $(0,\infty)\setminus (\alpha,\beta)$. (See again \cite[Lemma 3.1]{p}.) A function $f$ satisfying  \eqref{eq:G-class} which decreases on some interval of the positive real line is said to belong to the class $\mathcal G$. 

Suppose that $f$ has the representation \eqref{eq:G-class} with $r>0$. Then $f\in \mathcal G$ if and only if $(\log f)'(u)<0$ where $u>0$ is the point at which $(\log f)''(u)=0$. 

If $(\alpha,\beta)$ is the interval of decrease of $f\in \mathcal G$ (with $r>0$) then $u$ belongs to $(\alpha,\beta)$. It also follows that $\log f$ attains its minimum value on $[\beta,\infty)$ at $\beta$.

\begin{lemma}
\label{lemma:conformal-onto}
The holomorhic function $\log f$ maps $[\beta,\infty)$ onto $[\log f(\beta),\infty)$, and the vertical line $\beta+iy$, $y\in [0,\infty)$ is mapped to a curve  where both the real and imaginary parts decrease to $-\infty$ as $y\to \infty$. 

Furthermore, $\log f$ maps the domain $\{\Re z>\beta\} \cap \mathbb C_+$
conformally onto the domain $\Omega=\Omega_1\cup \Omega_2$, where
$$
\Omega_1=\{w\in \mathbb C\,|\, \Re w\geq \log f(\beta), \Im w>0\}
$$
and 
$$
\Omega_2=\{w\in \mathbb C\,|\, \exists y>0: \Re w=\Re \log f(\beta+iy), \Im w>\Im \log f(\beta+iy) \}.
$$
\end{lemma}
In the proof of this lemma we need a lower bound on $|\log f|$ for $f$ satisfying \eqref{eq:G-class} on quarter circles:
\begin{lemma}
\label{lemma:quarter-circle}Suppose that $f$ satisfies \eqref{eq:G-class}. Then 
$$
\inf\{ |\log f(Re^{i\theta})| \, | \, \theta\in[0,\pi/2]\} \to \infty \quad \text{as}\ R\to \infty.
$$
\end{lemma}
{\it Proof.} We denote by $n(t)=\#\{ k| \lambda_k\leq t\}$ the zero counting function associated with the sequence $\{\lambda_k\}$. Since $\sum_k1/\lambda_k^3<\infty$ and $\sum_k1/\lambda_k^2=\infty$ it follows that 
$$
\int_0^{\infty}\frac{n(t)}{t^4}\, dt<\infty, \quad \int_0^{\infty}\frac{n(t)}{t^3}\, dt=\infty \quad \text{and}\quad
\frac{n(t)}{t^3}\to 0 \text{ for } t\to \infty.
$$
Now consider the canonical product $P$ occuring in \eqref{eq:G-class}. For $z$ in the cut plane we can write $\log P(z)$ as a Stieltjes integral w.r.t.\ the zero counting function, and use integration by parts in order to obtain
$$
\log P(z)=z^3\int_0^{\infty}\frac{n(t)}{t^3}\frac{1}{t+z}\, dt.
$$
(See e.g.\ \cite[Lecture 12]{levin}.) This gives
$$
\log f(z)=r\log z+bz+z^2\left(\int_0^{\infty}\frac{n(t)}{t^3}\frac{z}{t+z}\, dt+a\right),
$$
so that 
$$
|\log f(z)|\geq |z^2|\left|\int_0^{\infty}\frac{n(t)}{t^3}\frac{z}{t+z}\, dt+a\right| -r|\log z|-|b||z|.
$$
The idea is to obtain a strictly positive lower bound on the factor containing the integral. A computation shows that (for $z=Re^{i\theta}$)
$$
\frac{z}{z+t}=\frac{R^2+Rt\cos \theta+iRt\sin \theta}{R^2+t^2+2Rt\cos \theta}.
$$
This yields (for $\theta\in [0,\pi/2]$) 
\begin{align*}
\left|\int_0^{\infty}\frac{n(t)}{t^3}\frac{z}{t+z}\, dt+a\right| &\geq
\int_0^{\infty}\frac{n(t)}{t^3}\Re \left(\frac{z}{t+z}\right)\, dt+a\\
&=\int_0^{\infty}\frac{n(t)}{t^3}\left(\frac{R^2+Rt\cos \theta}{R^2+t^2+2Rt\cos \theta}\right)\, dt+a\\
&\geq \int_0^{\infty}\frac{n(t)}{t^3}\left(\frac{R^2}{R^2+t^2+2Rt}\right)\, dt+a\\
&= \int_0^{\infty}\frac{n(t)}{t^3}\left(\frac{R}{R+t}\right)^2\, dt+a.
\end{align*}
By monotone convergence,
$$
\int_0^{\infty}\frac{n(t)}{t^3}\left(\frac{R}{R+t}\right)^2\, dt\uparrow \int_0^{\infty}\frac{n(t)}{t^3}\, dt=\infty.
$$
Since $|z^2|$ dominates both $|z|$ and $|\log z|$ the proof is completed.\hfill $\square$

{\it Proof of Lemma \ref{lemma:conformal-onto}.}
Let $\iota$ be the simple closed curve consisting of a straight line from $\beta$ to $R$, a circular part from $R$ to the point of intersection of the circle $|z|=R$ with the vertical line $x=\beta$ in the upper half plane, and finally the vertical line segment from the point of intersection to the point $\beta$.

The function $\log f$ is one-to-one in $\{\Re z\geq \beta\}\cap \{\Im z\geq0\}$: This can be seen using Lemma \ref{lemma:conformal} in $\{\Re z> u\}\cap \{\Im z>0\}$ and noting that \eqref{eq:logf-int} holds for $z_1,z_2\in \{\Re z\geq \beta\}\cap \{\Im z\geq0\}$ and $\log f$ is one-to-one in $[\beta,R]$. 
The image curve $\kappa$ of $\iota$ under the mapping $\log f$ is thus another simple closed curve. 

Let $w$ belong to the bounded region bounded by the curve $\kappa$. Then 
$$
\frac{1}{2\pi i}\int_{\iota}\frac{(\log f)'(z)}{\log f(z)-w}dz=\frac{1}{2\pi i}\int_{\kappa}\frac{1}{\zeta-w}d\zeta=1,
$$
and therefore the function $\log f-w$ has exactly one zero inside $\iota$. (This is a version of a theorem of Darboux, see e.g.\ \cite{osgood}.) Thus, the conformal image of the region bounded by $\iota$ under $\log f$ is the region bounded by $\kappa$. 

The image of the interval $[\beta,\infty)$ under $\log f$ is the interval $[\log f(\beta),\infty)$.
Let $\log f=U+iV$. Since $\partial_xU=\partial_yV$ and $U$ is increasing at $x=R$ we must have $V(R+iy)>0$ for $y>0$ close to $0$. Thus the circular arc of $\iota$ close to $R$ is mapped into the upper half plane. Next we consider the image of the vertical line $\beta+i(0,\infty)$. For $x>0$ we have from \eqref{eq:logf''}
\begin{align*}
\partial^2_{yy}V(x+iy)&=-\partial^2_{xx}V(x+iy)\\
&=-y\left(\frac{2xr}{(x^2+y^2)^2}+\sum_{k=1}^{\infty}\frac{2(x+\lambda_k)}{((x+\lambda_k)^2+y^2)^2}\right)<0.
\end{align*}
In particular, $\partial_yV(\beta+iy)$ is decreasing as a function of $y>0$ and thus 
$$
\partial_yV(\beta+iy)<\partial_yV(\beta)=\partial_xU(\beta)=(\log f)'(\beta)=0.
$$
This means that $V(\beta+iy)$ decreases from its value $0$ at $y=0$. We conlude that the curve $\log f(\beta+i(0,\infty))$ lies in the lower half plane. 

The relation \eqref{eq:imlogf'} yields
$$
\partial_{y}U(x+iy)=y\left(\frac{r}{x^2+y^2}-2a+\sum_{k=1}^{\infty}\left(\frac{1}{(x+\lambda_k)^2+y^2}-\frac{1}{\lambda_k^2}\right)\right),
$$
and so in particular $\partial_{y}U(\beta+iy)/y$ is decreasing as a function of $y>0$. Its value at $y=0$ is equal to 
$$
\frac{r}{\beta^2}-2a+\sum_{k=1}^{\infty}\left(\frac{1}{(\beta+\lambda_k)^2}-\frac{1}{\lambda_k^2}\right),
$$
which is exactly $-(\log f)''(\beta)$, a negative quantity because $\beta>u$. Therefore $\partial_{y}U(\beta+iy)<0$ for all $y>0$ and hence $U(\beta+iy)$ decreases for $y>0$. Both the real- and imaginary part of the curve $\log f(\beta+iy)$, $y>0$ are therefore decreasing functions. 

To finish the proof of the asserted description of the conformal image of $\{\Re z>\beta\}\cap \mathbb C_+$ it suffices to take an arbitrary $w\in \Omega$ and show that $w$ is inside the region bounded by $\kappa$ for sufficiently large $R$. By Lemma \ref{lemma:quarter-circle}, $|\log f(Re^{i\theta})|>|w|$ uniformly for $0\leq \theta\leq \pi/2$, when $R$ is chosen sufficiently large. 
The image of the circular part of $\iota$ is thus a curve in $\{|z|>|w|\}$ conneting the point $\log f(R)$ and a point in the lower half plane. As mentioned above the curve $\log f(\beta+iy)$ belongs to the upper half plane for small $y>0$ and by conformality it cannot cross the real half line $(\log f(\beta), \infty)$. This means that the curve must surround $w$. \hfill $\square$

\begin{thm}
\label{thm:G-main}
Let $f\in \mathcal G$ with $r>0$ and suppose that it has a local minimum at $\beta$, $\beta>0$. Then $f$ has an inverse function $f^{-1}$ defined on $(f(\beta),\infty)$ which can be extended to a univalent Pick-function defined in the cut plane $\mathbb C\setminus (-\infty,f(\beta)]$.  
\end{thm}
{\it Proof.} From Lemma \ref{lemma:conformal-onto} we know that $\log f$ is a conformal mapping of $\{\Re z>\beta\}\cap \mathbb C_+$ onto $\Omega$. Since the strip $S=\{0<\Im w<\pi\}$ is a subset of $\Omega$ it follows that $f^{-1}\equiv (\log f)^{-1}(\log w)$ maps the upper half plane into the domain $\{\Re z>\beta\}\cap \mathbb C_+$, and in particular into the upper half plane. Hence $f^{-1}$ is a Pick-function and we also have  $f((\log f)^{-1}(\log w))=w$. \hfill $\square$

\begin{rem}
Not only $f^{-1}$ but also the function $w\mapsto (f^{-1}(w)-\beta)^2$ is a Pick function. This is seen noting that $f^{-1}$ maps the upper half plane into the quardrant $\{ \Re z>\beta, \Im z>0\}$.
\end{rem}

\begin{rem}
If $(\log f)'(u)=0$ Theorem \ref{thm:G-main} still holds (with $\beta=u$). However, if $(\log f)'(u)>0$ the conformal image of the quardrant $\{ \Re z>u\}\cap \mathbb C_+$ does not cover the entire strip $\{ 0<\Im w<\pi\}$ and in this situation it is not possible to extend the inverse to a Pick-function.
\end{rem}

\begin{prop}
Let $f\in \mathcal G$ with $r>0$ and let $\beta>0$ be a local minimum point for $f$. Let $f^{-1}$ be the inverse defined on $(f(\beta), \infty)$ and extended to a Pick function. Then, for $w\in \mathbb C\setminus (-\infty, f(\beta)]$,
$$
f^{-1}(w)=\beta+(w-f(\beta))\int_{0}^{\infty}\frac{d(t)}{t+w-f(\beta)}\, dt-\frac{c}{w}
$$
where $c\geq 0$, and $d(s)$ is positive and real analytic on $(0,f(\beta))\cup (f(\beta),\infty)$ satisfying $\int_0^{\infty}d(t)/(t+1)<\infty$.
\end{prop}
{\it Sketch of Proof.} Since $w\mapsto f^{-1}(w+f(\beta))-\beta$ is a Pick-function which can be analytically extended across the positive real line, and which maps the positive real line into itself it has the representation  
$$
f^{-1}(w+f(\beta))-\beta=aw+w\int_{0}^{\infty}\frac{d\sigma(t)}{t+w}
$$
where $a=\lim_{x\to \infty}f^{-1}(x+f(\beta))/x\geq 0$, and $\sigma$ is a positive measure on $[0,\infty)$ satisfying $\int_0^{\infty}d\sigma(t)/(t+1)<\infty$. See \cite[Theorem 2.4]{b}. It readily follows that $a=0$. It can also be shown that $\sigma$ has an analytic density on $(0,f(\beta))\cup (f(\beta),\infty)$, and that this leaves the possibility of a point mass only at $f(\beta)$.
\hfill $\square$

As also mentioned in \cite{p}, the motivation for studying the class $\mathcal G$ and the results above comes from the multiple gamma functions introduced by Barnes, see \cite{b1}. The double gamma function is usually denoted by $\Gamma_2$ and it is defined as $\Gamma_2(z)=(2\pi)^{z/2}/G(z)$, where 
$$
G(z+1)=(2\pi)^{z/2}e^{-((1+\gamma)z^2+z)/2}\prod_{k=1}^{\infty}\left(
  1+\frac{z}{k}\right) ^ke^{-z+z^2/2k}.
$$
The entire function $G$ is the $G$-function of Barnes and it satisfies the relation $G(z+1)=\Gamma(z)G(z)$. 

\begin{cor}  The following hold
\begin{enumerate}
\item[(a)] Barnes $G$-function belongs to $\mathcal G$ and increases on the interval $(\beta_G,\infty)$. Its inverse maps $(G(\beta_G),\infty)$ to $(\beta_G,\infty)$ and it can be extended to a univalent Pick-function. 
\item[(b)] The reciprocal to Barnes double gamma function $\Gamma_2$ belongs to $\mathcal G$ and increases on the interval $(\beta_2,\infty)$. Its inverse maps $(1/\Gamma_2(\beta_2),\infty)$ to $(\beta_2,\infty)$ and it can be extended to a univalent Pick-function. 
\end{enumerate}
\end{cor} 
According to \cite{o} and \cite{p} both $G$ and $1/\Gamma_2$ belong to the class $\mathcal G$ and therefore the corollary follows.
Numerical experiments yield that 
$\beta_G\approx 2.568$, $G(\beta_G)\approx 0.945$, 
$\beta_2\approx 3.763$, and $1/\Gamma_2(\beta_2)\approx 0.048$. 

We end this paper with a result concerning the inverse of the double gamma function itself. The proof is almost the same as the proof in \cite[Proposition 2.6]{p}, so we leave out the details.
\begin{prop}Barnes' double gamma function $\Gamma_2$ decreases on $(\beta_2,\infty)$ and maps it onto $(0,\Gamma_2(\beta_2))$. Its inverse $k$ is of the form $k(w)=h(1/w)$, $h$ being a Pick-function with an  extension across $(1/\Gamma_2(\beta_2),\infty)$. Furthermore, $k$ can be represented as 
$$
k(w)=k(\Gamma_2(\beta_2/2))-\int_{0}^{\Gamma_2(\beta)}\frac{w-\Gamma_2(\beta)/2}{w+s-2ws/\Gamma_2(\beta)}\,
d\nu(s),
$$
where $\nu$ is a positive measure on $[0,\Gamma_2(\beta_2)]$.
\end{prop}

\noindent Department of Mathematics\\University of Copenhagen\\Universitetsparken 5, 2100 Copenhagen, Denmark.

\end{document}